\documentclass{article}

\usepackage{arxiv}

\usepackage[utf8]{inputenc} 
\usepackage[T1]{fontenc}    
\usepackage{hyperref}       
\usepackage{url}            
\usepackage{booktabs}       
\usepackage{amsfonts}       
\usepackage{nicefrac}       
\usepackage{microtype}      
\usepackage{lipsum}
\usepackage{graphicx}
\graphicspath{ {./images/} }
\usepackage{amsmath}
\usepackage{float}
\usepackage{url}
\newtheorem{theorem}{Theorem}[section]      
\newtheorem{lemma}[theorem]{Lemma}          

\title{A New Special Function and Its Application in Probability
}

\author{
 Zeraoulia Rafik \\
  University of Khemis Miliana,Algeria\\
Deprtaement of mathematics\\
 Laboratory of Pure and Applied Mathematics (LMPA)  \\
  \texttt{zeraoulia@univ-dbkm.dz} \\
   \And
 Alvaro H Salas\\
  Universidad Nacional de Colombia,Colombia\\
  Departemnt of physisc\\
 Fizmako Group research\\
  \texttt{ahsalass@unal.edu.co} \\
  \And
 David L. Ocampo \\
 Universidad Nacional de Colombia,Colombia\\
  Universidad Caldas,Colombia\\
  departement of mathematics \\
  \texttt{david.ocampo@ucaldas.edu.co} \\
}

\begin{document}
\maketitle
\begin{abstract}
 In this note we present a new special function that behaves like the error function and we provide an approximated accurate closed form for its CDF in terms of both Chebyshev polynomials of the first kind and the error function. Also, we provide its series representation using Pad\'e approximant. We show convincing numerical evidence of an accuracy of $10^{-6}$ for the approximants in the sense of the quadratic mean norm. A similar approach may be applied to other probability distributions, for example, the Maxwell--Boltzmann distribution and the normal distribution, and we show its application using both of those distributions.
\end{abstract}


\section{Introduction}

Integrals of the error function, see (1), occur in a great variety of applications usually in problems involving multiple integration where the integrand contains exponentials of the squares of the argument; an example of applications can be cited from atomic physics astrophysics and statistical analysis. It comes into our mind to seek for the integration of such functions $f(x)$ power its antiderivative $g(x)$. We have got example (1) where it is the power of two distributions related to normal distribution \cite{TeugelsSundt2004} as shown below such that $f(x)=e^{-x^{2}}$ and $g(x)=\operatorname{erf}(x)$:
\begin{equation}
I(a)=\int_{0}^{a} \left(e^{-x^{2}}\right)\operatorname{erf}(x)\,dx.
\tag{1}
\end{equation}
The function $\operatorname{erf}(x)$ is called the error function and it is defined in (2):
\begin{equation}
\frac{2}{\sqrt{\pi}}\int_{0}^{x} e^{-t^{2}}\,dt=\operatorname{erf}(x).
\tag{2}
\end{equation}

\subsection{Numerical approximation of $\displaystyle\int_{0}^{a} \left(e^{-x^{2}}\right)\operatorname{erf}(x)\,dx$ in some ranges of values}

Now, if we really need a simple expression for $I(a)$ in some range of values, there are ways to get various approximations.

The function is very nice. It goes to its limit at $\infty$ very fast.

\begin{figure}[H]
  \centering
  \includegraphics[width=0.6\textwidth]{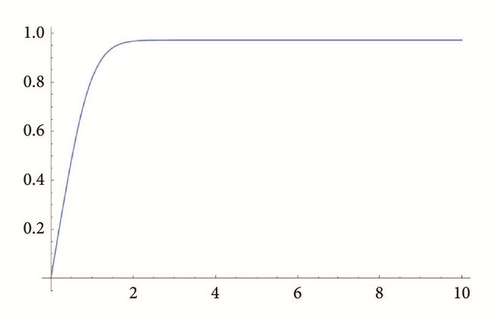}
  \caption{The plot of $I(a)$ for $a \in [0,10]$.}
\end{figure}

Therefore (depending on the accuracy we need) we can easily take $I(a)=I(\infty)$ for $a>a_{0}$ with $a_{0}$ around $3$ or $4$. Mathematica gives the following for the first 100 digits:
\begin{equation}
\begin{split}
I(\infty) &= 0.972106992769178593151077\\
&\quad 875442391175554272183385569900\\
&\quad 9722910408441888759958220033410678218401258734 \cdots
\end{split}
\tag{3}
\end{equation}

Now, what can we do for small $a$?

The function is so nice; we can just use the Taylor expansion around $a=0$. The first term is as follows:
\begin{equation}
I(a)\approx a.
\tag{4}
\end{equation}

\begin{figure}[H]
  \centering
  \includegraphics[width=0.6\textwidth]{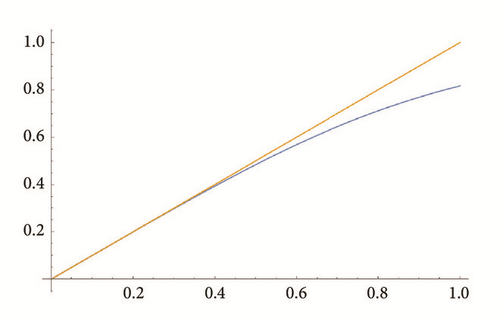}
  \caption{Approximation of $I(a)$ for $a \in [0,1]$ using Taylor expansion.}
\end{figure}

The proof is simple. The Taylor series looks like the following:
\begin{equation}
I(a)=I(0)+I'(0)a+\frac{I''(0)}{2!}a^{2}+\frac{I^{(3)}(0)}{3!}a^{3}+\cdots
\tag{5}
\end{equation}
We may see the following:
\begin{equation}
I(0)=0,\qquad I'(0)=e^{-a^{2}}\operatorname{erf}(a)\big|_{a=0}=1.
\tag{6}
\end{equation}
Now let us find a better approximation by computing the higher derivatives:
\begin{equation}
I''(a)=\left(e^{-a^{2}}\operatorname{erf}(a)\right)'
=-\frac{2}{\sqrt{\pi}}\,a e^{-a^{2}}\left(\operatorname{erf}(a)+1\right)\left(\sqrt{\pi} e^{a^{2}}\operatorname{erf}(a)+a\right),
\qquad I''(0)=0.
\tag{7}
\end{equation}
We use Mathematica as a shortcut, but it is easy to do it by hand, if we remember that
\[
\operatorname{erf}'(x)=\frac{2}{\sqrt{\pi}}e^{-x^{2}}.
\]
We obtain
\begin{equation}
I^{(3)}(0)=0,\qquad I^{(4)}(0)=-\frac{12}{\sqrt{\pi}}.
\tag{8}
\end{equation}
So our next approximation is as follows:
\begin{equation}
I(a)\approx a-\frac{1}{2\sqrt{\pi}}a^{4}.
\tag{9}
\end{equation}

\begin{figure}[H]
  \centering
  \includegraphics[width=0.6\textwidth]{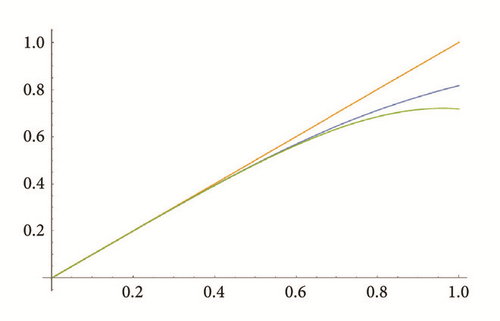}
  \caption{Plot of $I(a)$ with both approximations.}
\end{figure}

Another way to approximate the function \cite{MaiDuyTranCong2003} is using its derivative:
\begin{equation}
\frac{dI}{da}=e^{-a^{2}}\operatorname{erf}(a).
\tag{10}
\end{equation}
But this is an ordinary differential equation, which can be solved numerically.

As an illustration, here is a simple explicit Euler scheme for the step size $h$:
\begin{equation}
\frac{I(a+h)-I(a)}{h}=e^{-a^{2}}\operatorname{erf}(a),
\qquad
I(a+h)=I(a)+h\,e^{-a^{2}}\operatorname{erf}(a).
\tag{11}
\end{equation}
We can use an initial value $I(0)=0$.

For $h=1/10$, we have the following result (red dots) compared to the exact function (blue line).

\begin{figure}[H]
  \centering
  \includegraphics[width=0.6\textwidth]{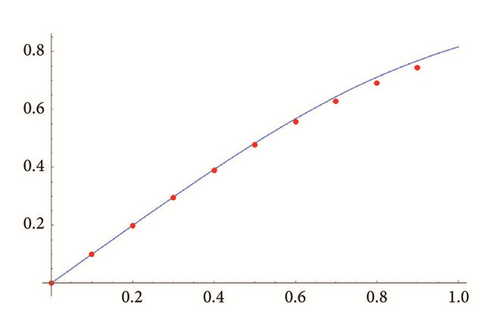}
  \caption{A simple explicit Euler scheme for the step size $h=1/10$.}
\end{figure}

For $h=1/50$ see the next figure.

\begin{figure}[H]
  \centering
  \includegraphics[width=0.6\textwidth]{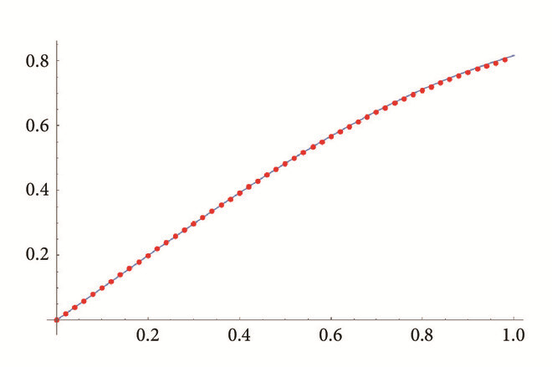}
  \caption{A simple explicit Euler scheme for the step size $h=1/50$.}
\end{figure}

This way can serve as a good alternative to numerical integration \cite{SuzukiSuzuki2003} (depending on the context and the application of course). Let us now show the relationship between this function and other standard special functions (integral of error function) \cite{NgGeller1969} as error function and cumulative distribution function for normal distribution in the context of its use. Function (1) could be used to find values of complicated integral which are not available in any references of standard special functions and also it is not available to get their values in Wolfram Alpha, for example,
\begin{equation}
\int_{0}^{+\infty} e^{x^{2}}\bigl(1-2\Phi(x\sqrt{2})\bigr)\,dx,
\tag{12}
\end{equation}
with $\Phi(x)=(1/\sqrt{2\pi})\int_{-\infty}^{x} e^{-z^{2}/2}\,dz$. It is the CDF (cumulative distribution function) for the normal distribution; if someone was asked to find the value of this integral, he would be confused because it is very complicated; probably he cannot show whether it is convergent or not; even Wolfram Alpha as a best means of computation cannot recognize at least that $\Phi$ is a cumulative normal distribution, so no result would be obtained about the value of this integral. Let us compute (12) using (1) and we will conclude that they have the same value and both are identical function and identical integral.

The well-known formula which expresses the relationship between the error function and the cumulative density function, see (2), is defined as
\begin{equation}
\operatorname{Erf}(x)=2\bigl(\Phi(x\sqrt{2})-\Phi(0)\bigr)
=2\Bigl(\Phi(x\sqrt{2})-\frac{1}{2}\Bigr)
=2\Phi(x\sqrt{2})-1.
\tag{13}
\end{equation}
And it is easy to check that it always holds for every real number by the following short proof.

\paragraph{Proof.}
By definition, the error function
\begin{equation}
\operatorname{Erf}(x)=\frac{2}{\sqrt{\pi}}\int_{0}^{x} e^{-t^{2}}\,dt.
\tag{14}
\end{equation}
Writing $t^{2}=z^{2}/2$ implies $t=z/\sqrt{2}$ (because $t$ is not negative), whence $dt=dz/\sqrt{2}$. The endpoints $t=0$ and $t=x$ become $z=0$ and $z=x\sqrt{2}$. To convert the resulting integral into something that looks like a cumulative distribution function (CDF), it must be expressed in terms of integrals that have lower limits of $-\infty$; thus
\begin{equation}
\operatorname{Erf}(x)
=\frac{2}{\sqrt{2\pi}}\int_{0}^{x\sqrt{2}} e^{-z^{2}/2}\,dz
=2\left(\frac{1}{\sqrt{2\pi}}\int_{-\infty}^{x\sqrt{2}} e^{-z^{2}/2}\,dz
-\frac{1}{\sqrt{2\pi}}\int_{-\infty}^{0} e^{-z^{2}/2}\,dz\right).
\tag{15}
\end{equation}
Those integrals on the right hand side are both values of the CDF of the standard normal distribution:
\begin{equation}
\Phi(x)=\frac{1}{\sqrt{2\pi}}\int_{-\infty}^{x} e^{-z^{2}/2}\,dz.
\tag{16}
\end{equation}
Specifically,
\begin{equation}
\operatorname{Erf}(x)=2\bigl(\Phi(x\sqrt{2})-\Phi(0)\bigr)
=2\Bigl(\Phi(x\sqrt{2})-\frac{1}{2}\Bigr)
=2\Phi(x\sqrt{2})-1.
\tag{17}
\end{equation}

Now since the LHS of (18) has a known value which is $0.97210699\cdots$, then the right hand side also equals $0.97210699\cdots$; hence we came up with the following identity:
\begin{equation}
\int_{0}^{a} \left(e^{-x^{2}}\right)\operatorname{Erf}(x)\,dx
=\int_{0}^{a} e^{x^{2}}\bigl(1-2\Phi(x\sqrt{2})\bigr)\,dx.
\tag{18}
\end{equation}

Now we shall call the function defined in (1)
\[
T(x)=\int_{0}^{x} \left(e^{-t^{2}}\right)\operatorname{erf}(t)\,dt
\]
since it does not refer to anyone and it has unknown analytic representation as elementary function using standard special functions and the RHS of (18) presents another representation of $T(x)$ function using CDF of the normal distribution.

\begin{lemma}
$T(x)=\displaystyle\int_{0}^{x} \left(e^{-t^{2}}\right)\operatorname{erf}(t)\,dt$ cannot be expressed in terms of elementary function.
\end{lemma}

\paragraph{Proof.}
It is a theorem of Liouville \cite{Liouville1838}, reproven later with purely algebraic methods, that for rational functions $f$ and $g$, $g$ nonconstant, the antiderivative
\begin{equation}
\int \bigl[f(x)\exp(g(x))\bigr]\,dx
\tag{19}
\end{equation}
can be expressed in terms of elementary functions if and only if there exists some rational function $h$ such that it is a solution to the differential equation:
\begin{equation}
f=h'+hg.
\tag{20}
\end{equation}
Now if we apply Liouville theorem we can come up with the following ODE:
\[
1=h'(x)+h(x)\bigl(-x^{2}\operatorname{erf}(x)\bigr)
\]
with $g(x)=-x^{2}\operatorname{erf}(x)$ and $f(x)=1$. It is first ordinary differential equation. The computation we made with Wolfram Alpha gives the following solution:
\begin{equation}
\begin{split}
&\exp\!\Biggl(
\frac{e^{-x^{2}}\Bigl(e^{x^{2}}(\sqrt{\pi}x^{3}\operatorname{erf}(x)-1)+x^{2}+1\Bigr)}{3\sqrt{\pi}}
\Biggr) \\
&\quad \times
\Biggl(
\exp\Biggl\{
\frac{1}{3}\sqrt{\pi}\Biggl(
\int_{1}^{x} 
\exp\!\Bigl(
\frac{1}{3}\Bigl(t^{3}(-\operatorname{erf}(t)) 
- \frac{e^{-t^{2}}(t^{2}+1)}{\sqrt{\pi}}\Bigr)
\Bigr)\,dt
- \int_{0}^{1} 
\exp\!\Bigl(
\frac{1}{3}\Bigl(t^{3}(-\operatorname{erf}(t)) 
- \frac{e^{-t^{2}}(t^{2}+1)}{\sqrt{\pi}}\Bigr)
\Bigr)\,dt
\Biggr)
\Biggr\} + 1
\Biggr)
\end{split}
\tag{21}
\end{equation}

with $h(0)=1$. Really the function $h$ can be written as follows:
\begin{equation}
h(x)=\ell(x)\left[c_{1}+\int_{x}^{1} \ell(-\xi)\,d\xi\right].
\tag{22}
\end{equation}
Now it is clear that $\ell(x)$ is a transcendental function and the defined integral in the right hand side of the $h(x)$ expression is also transcendental function because we have derivatives of rational functions being rational functions. Therefore, if the antiderivative is rational, then the original function was rational. The function $h$ is rational only at $x=0$, and since $h(x)\neq 0$, then the sum of two transcendental functions is always transcendental function. According to definition of the rational function, $h(x)$ cannot be called a rational function; then we are done.

\section{A Possible Approach Formula for $T(+\infty)$}

We may give here a possible approach formula for $T(+\infty)$ which is defined as follows:
\begin{equation}
T(+\infty)=\int_{0}^{+\infty} \exp\bigl(-x^{2}\operatorname{erf}(x)\bigr)\,dx
=0.97210699\cdots.
\tag{23}
\end{equation}
The inverse symbolic calculator is unable to give us the representation of $0.97210699\cdots$ using standard special functions, but we have tried to give its representation using error function representation as hypergeometric function \cite{arxiv1702.08438}; we have
\begin{equation}
\operatorname{erf}(x)=\frac{2}{\sqrt{\pi}}\,x\,{}_1F_{1}\!\left(\frac{1}{2};\frac{3}{2};-x^{2}\right)
\tag{24}
\end{equation}
with ${}_1F_{1}$ being the Kummer confluent hypergeometric function \cite{arxiv1702.08438}. Now we have from (24) the following:
\begin{equation}
\int_{0}^{+\infty} \exp\bigl(-x^{2}\operatorname{erf}(x)\bigr)\,dx
=\int_{0}^{+\infty} \exp\!\left(-\frac{2}{\sqrt{\pi}}x^{3}\right)\,
{}_1F_{1}\!\left(\frac{1}{2};\frac{3}{2};-x^{2}\right)\,dx.
\tag{25}
\end{equation}
The RHS of (25) using (24) gives
\[
\bigl(\pi^{1/6}\Gamma(4/3)/2^{1/3}\bigr)\,{}_1F_{1}\!\left(0.5;1.5;-x^{2}\right).
\]
Hence we may choose $x=e\sqrt{\pi}$ and we can get finally the following:
\begin{equation}
T(+\infty)\sim
\frac{\pi^{1/6}\Gamma(4/3)}{2^{1/3}}\,
{}_1F_{1}\!\left(0.5;1.5;-\pi e^{2}\right)
=0.97216864\cdots.
\tag{26}
\end{equation}

\begin{figure}[H]
  \centering
  \includegraphics[width=0.6\textwidth]{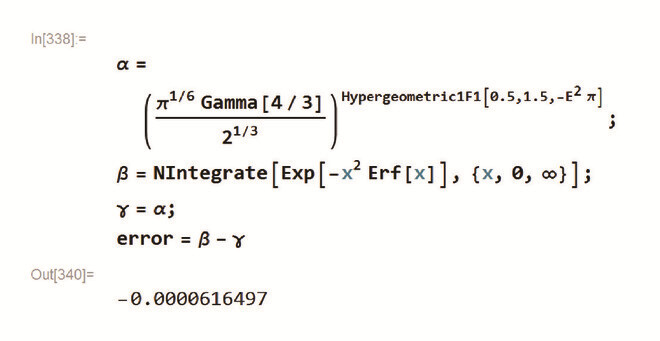}
  \caption{Error approximation for $T(+\infty)$.}
\end{figure}

Mathematica gives the nice approximation of (25) as shown in Figure~6.

\section{Series Representation of $T(x)$ Function}

We may try to find a series expansion in powers of $t$ of
\begin{equation}
I(t)=\int_{0}^{t} \exp\bigl(-x^{2}\operatorname{erf}(x)\bigr)\,dx
=\sum_{p=1}^{\infty} c_{p} t^{p}.
\tag{27}
\end{equation}
The coefficients $c_{p}
= d_{p-1}/p$ follow from the series expansion
\[
e^{-x^{2}}\operatorname{erf}x=\sum_{p=0}^{\infty} d_{p}x^{p},
\]
resulting in
\begin{equation}
I(t)=\sum_{p=1}^{\infty} c_{p} t^{p}
=t-\frac{t^{4}}{2\sqrt{\pi}}
+\frac{t^{6}}{9\sqrt{\pi}}
+\frac{2t^{7}}{7\pi}
-\frac{t^{8}}{40\sqrt{\pi}}
-\frac{4t^{9}}{27\pi}
+\frac{(\pi-28)t^{10}}{210\pi^{3/2}}
+O(t^{11}).
\tag{28}
\end{equation}
The series $I(1)=\sum_{p=1}^{\infty} c_{p}$ seems to converge. See Figure~7.

\begin{figure}[H]
  \centering
  \includegraphics[width=0.6\textwidth]{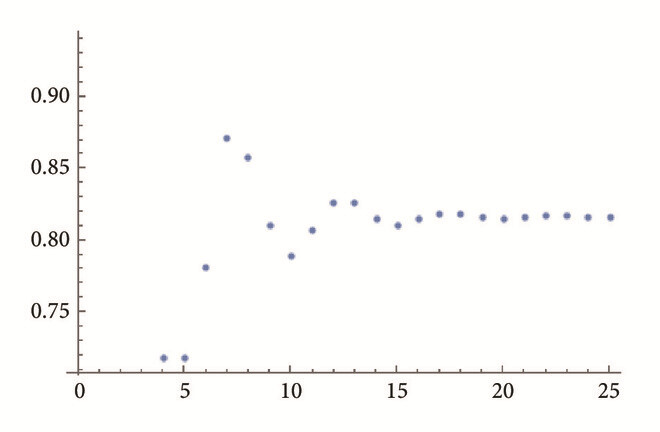}
  \caption{Convergence of $I_{N}=\sum_{p=0}^{N} c_{p}$ as a function of $N$ up to $N=25$.}
\end{figure}

The value of $I_{25}=0.8162$ agrees with $I(1)=0.816377$ to three decimal places. For $N=50$ the agreement is up to six decimal places, but this did not give us the power series closed form for $n$th term. We should use some approximations using approximation of the error function and Pad\'e approximant as shown in the following sections.

\section{Series Expansion of the $n$-th Derivative of $T(x)
=\displaystyle\int_{0}^{x} e^{-\xi^{2}}\operatorname{erf}(\xi)\,d\xi$}

\begin{lemma}
Series expansion of $T(x)=\displaystyle\int_{0}^{x} e^{-\xi^{2}}\operatorname{erf}(\xi)\,d\xi$ is
defined by this identity:
\begin{equation}
\int_{0}^{x} e^{-\xi^{2}}\operatorname{erf}(\xi)\,d\xi
=\sum_{n=0}^{\infty}
\lim_{\varepsilon\to 0}
\left(
\sum_{\substack{k_{1}+2k_{2}+\cdots+nk_{n}=n\\ k_{1}\ge 0,k_{2}\ge 0,\dots,k_{n}\ge 0}}
\prod_{j=1}^{n}
\frac{A_{j,\varepsilon}^{\,k_{j}}}{k_{j}!}
\right)\frac{x^{n+1}}{n+1},
\tag{29}
\end{equation}
where
\begin{equation}
A_{j,\varepsilon}=
\begin{cases}
\dfrac{2(-1)^{(j-1)/2}}{(j-2)\bigl(\frac{1}{2}(j-3)\bigr)!\sqrt{\pi}},
& \text{if } j\ge 3 \text{ and $j$ an odd integer},\\[1ex]
\varepsilon, & \text{otherwise }(0<\varepsilon<1),
\end{cases}
\tag{30}
\end{equation}
which is the key idea to get.
\end{lemma}

\paragraph{Proof.}
Suppose that we have the Taylor expansions:
\begin{equation}
f(x)=\sum_{n=1}^{\infty}\frac{a_{n}}{n!}x^{n}
\tag{31}
\end{equation}
and
\begin{equation}
g(x)=\sum_{n=1}^{\infty}\frac{b_{n}}{n!}x^{n}.
\tag{32}
\end{equation}
Then we have the standard result:
\begin{equation}
g(f(x))=\sum_{n=1}^{\infty}
\left(
\sum_{k=1}^{n} b_{k} B_{n,k}(a_{1},\dots,a_{n-k+1})
\right)\frac{x^{n}}{n!},
\tag{33}
\end{equation}
where $B_{n,k}(\cdot)$ are the partial Bell polynomials, which are defined by the following formula:
\begin{equation}
\widehat{B}_{m,j}(x_{1},x_{2},\dots,x_{m-j+1})
=\sum_{\substack{k_{0}+k_{1}+\cdots+k_{N}=j\\ k_{1}+2k_{2}+\cdots+Nk_{N}=m}}
\binom{j}{k_{0},k_{1},\dots,k_{N}}
\prod_{i=1}^{N} x_{i}^{k_{i}}.
\tag{34}
\end{equation}
The key idea to get series expansion of the $n$-th derivative of $T(x)$ which is defined in (29) is to use Taylor expansion of $g(x)=\exp(-x)$ and $f(x)=-x^{2}\operatorname{erf}(x)$ coming up for using one of the important formulas in mathematics called Bruno–Fadi formula such as that defined above in (33) using (34). It is well known that the Taylor expansion of $\exp(-x)$ is given by the following:
\begin{equation}
g(x)=\exp(-x)=\sum_{n=0}^{\infty}\frac{(-1)^{n}}{n!}x^{n}.
\tag{35}
\end{equation}

Probably the interesting here for readers to know is Taylor expansion of $\operatorname{erf}(x)$; we give a simple proof about its expansion series using Hermite polynomial.

\begin{lemma}
Taylor expansion of $\operatorname{erf}(x)$ at each point $a$ is given by this identity:
\begin{equation}
\operatorname{erf}_{a}(x)
=e^{-a^{2}}\sum_{n=0}^{\infty}\frac{(-1)^{n}H_{n}(a)}{n!}(x-a)^{n},
\tag{36}
\end{equation}
with $H_{n}(a)$ the Hermite polynomial of degree $n$.
\end{lemma}

\paragraph{Proof.}
$f^{(n)}(a)$ can be written in terms of Hermite polynomials $H_{n}$:
\begin{equation}
H_{0}(x)=1,\quad
H_{1}(x)=2x,\quad
H_{2}(x)=4x^{2}-2,\quad
H_{3}(x)=8x^{3}-12x,\quad
H_{4}(x)=16x^{4}-48x^{2}+12,\quad
H_{5}(x)=\cdots.
\tag{37}
\end{equation}
We may recognize that $H_{2n-1}(0)=0$, which gives the power series for $e^{-x^{2}}$ at $a=0$:
\begin{equation}
e^{-x^{2}}=1-\frac{2}{2!}x^{2}+\frac{12}{4!}x^{4}-\frac{120}{6!}x^{6}+\cdots.
\tag{38}
\end{equation}
After multiplying by $2/\sqrt{\pi}$, this integrates to
\begin{equation}
\operatorname{erf}(z)
=\frac{2}{\sqrt{\pi}}\left(z-\frac{z^{3}}{3}+\frac{z^{5}}{10}-\frac{z^{7}}{42}+\frac{z^{9}}{216}-\cdots\right).
\tag{39}
\end{equation}
Since $\dfrac{d^{n}}{dx^{n}}e^{-x^{2}}=(-1)^{n}e^{-x^{2}}H_{n}(x)$, one can do a Taylor series for every $a$:
\begin{equation}
\operatorname{erf}_{a}(x)
=e^{-a^{2}}\sum_{n=0}^{\infty}\frac{(-1)^{n}H_{n}(a)}{n!}(x-a)^{n}.
\tag{40}
\end{equation}
Then we are done.
Now by composition of (40) with (35) after multiplying (40) by the term $-x^{2}$, we come up to Bruno–Fadi formula which is defined as
\begin{equation}
e^{-x^{2}}\operatorname{erf}(x)
=\sum_{n=0}^{\infty}
\lim_{\varepsilon\to 0}
\left(
\sum_{\substack{k_{1}+2k_{2}+\cdots+nk_{n}=n\\ k_{1}\ge 0,k_{2}\ge 0,\dots,k_{n}\ge 0}}
\prod_{j=1}^{n}\frac{A_{j,\varepsilon}^{\,k_{j}}}{k_{j}!}
\right)x^{n},
\tag{41}
\end{equation}
where
\begin{equation}
A_{j,\varepsilon}=
\begin{cases}
\dfrac{2(-1)^{(j-1)/2}}{(j-2)\bigl(\frac{1}{2}(j-3)\bigr)!\sqrt{\pi}},
& \text{if } j\ge 3 \text{ and $j$ an odd integer},\\[1ex]
\varepsilon, & \text{otherwise }(0<\varepsilon<1).
\end{cases}
\tag{42–43}
\end{equation}
Integrating this equation term by term gives
\begin{equation}
\int_{0}^{x} e^{-\xi^{2}}\operatorname{erf}(\xi)\,d\xi
=\sum_{n=0}^{\infty}
\lim_{\varepsilon\to 0}
\left(
\sum_{\substack{k_{1}+2k_{2}+\cdots+nk_{n}=n\\ k_{1}\ge 0,k_{2}\ge 0,\dots,k_{n}\ge 0}}
\prod_{j=1}^{n}\frac{A_{j,\varepsilon}^{\,k_{j}}}{k_{j}!}
\right)\frac{x^{n+1}}{n+1},
\tag{44}
\end{equation}
which gives the series expansion for the new special function.

Using Mathematica as a shortcut, it shows that (44) holds and also it shows the incrementation of $\pi$ as shown in Figure~8. We may add also the series expansion of $T(x)$, the new special function, using Mathematica code as shown in Figure~9.

\begin{figure}[H]
  \centering
  \includegraphics[width=0.6\textwidth]{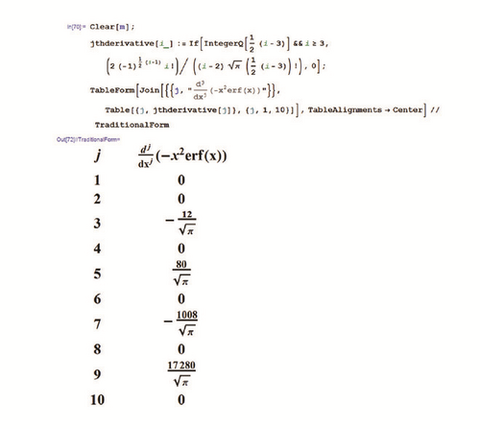}
  \caption{$n$-th derivative of new special function $T(x)$.}
\end{figure}

\begin{figure}[H]
  \centering
  \includegraphics[width=0.6\textwidth]{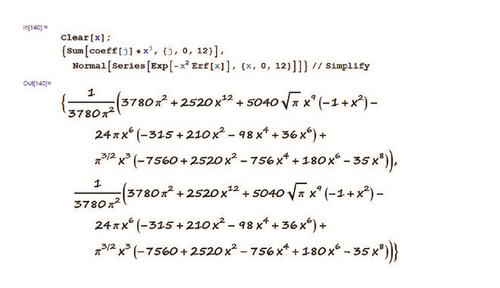}
  \caption{Series expansion of new special function $T(x)$ around $x=0$.}
\end{figure}

\section{Series Representation of $\displaystyle\int_{-1}^{1} \operatorname{erf}(x)^{n}\,dx$ Using Error Function Approximation}

We have the power series of
\begin{equation}
e^{-x^{2}}\operatorname{erf}(x)
=\sum_{k=0}^{\infty}(-1)^{k}\frac{x^{2k}\operatorname{erf}^{(k)}(x)}{k!},
\tag{45}
\end{equation}
and then from (45) we have the following:
\begin{equation}
\int_{-1}^{1} e^{-x^{2}}\operatorname{erf}(x)\,dx
=\sum_{k=0}^{\infty}\frac{(-1)^{k}}{k!}
\int_{-1}^{1} x^{2k}\operatorname{erf}^{(k)}(x)\,dx.
\tag{46}
\end{equation}
Now it is hard so much to evaluate the integral in RHS of (46) using error function expression; then we should use the following nice approximation:
\begin{equation}
\bigl(\operatorname{erf}(x)\bigr)^{2}\approx 1-e^{-a x^{2}}
\quad\text{with}\quad
a=\frac{(1+\pi)^{2/3}}{\log^{2}(2)}.
\tag{47}
\end{equation}
Here we can give a short proof to show that the error function squared was approximated as well with the value of $a=(1+\pi)^{2/3}\log^{2}(2)$.

\paragraph{Proof.}
We fully agree that
\begin{equation}
F(a)=\int_{0}^{\infty}
\bigl(\operatorname{erf}(x)^{2}-(1-e^{-a x^{2}})\bigr)^{2}\,dx
\tag{48}
\end{equation}
is minimum for $a\approx 1.23907$. According to RIES, this number seems to be much closer to
\begin{equation}
a=(1+\pi)^{2/3}\log^{2}(2)\approx 1.23907
\tag{49}
\end{equation}
than to $\pi^{2}/8\approx 1.23370$ even if this does not make very large difference (the maximum error is reduced from $0.006$ to $0.004$ and the value of the integral $F(a)$ changes from $0.00002769$ to $0.00002572$). If we look for a still better approximation, we could consider $\log(1-\operatorname{erf}(x)^{2})$ (which, for sure, introduces a bias in the problem), establish a Pad\'e approximant, and finally arrive to
\begin{equation}
\operatorname{erf}(x)^{2}\approx
1-\exp\!\left(-\frac{4}{\pi}\,\frac{1+\alpha x^{2}}{1+\beta x^{2}}\,x^{2}\right)
\tag{50}
\end{equation}
where
\begin{equation}
\alpha=\frac{10-\pi^{2}}{5(\pi-3)\pi},
\qquad
\beta=\frac{120-60\pi+7\pi^{2}}{15(\pi-3)\pi}.
\tag{51}
\end{equation}
The value of the corresponding error function is $1.1568\times 10^{-7}$, that is to say, almost $250$ times smaller than that with the initial formulation; the maximum error is $0.00035$.

Now we are ready to approximate
\begin{equation}
I_{n}=\int_{-1}^{1}\bigl(\operatorname{erf}(x)\bigr)^{2n}\,dx
\tag{52}
\end{equation}
by
\begin{equation}
J_{n}=\int_{-1}^{1}\bigl(1-e^{-a x^{2}}\bigr)^{n}\,dx,
\tag{53}
\end{equation}
for which the binomial expansion would be required (easy). This would give you things like the following:
\begin{align}
J_{1}&=2-\frac{\sqrt{\pi}}{\sqrt{a}}\operatorname{erf}(\sqrt{a}),
\tag{54}\\[0.5ex]
J_{2}&=2-\frac{2\sqrt{\pi}}{\sqrt{a}}\operatorname{erf}(\sqrt{a})
+\frac{\sqrt{\pi/2}}{\sqrt{a}}\operatorname{erf}(\sqrt{2a}),
\tag{55}\\[0.5ex]
J_{3}&=2-\frac{3\sqrt{\pi}}{\sqrt{a}}\operatorname{erf}(\sqrt{a})
+\frac{3\sqrt{\pi/2}}{\sqrt{a}}\operatorname{erf}(\sqrt{2a})
-\frac{\sqrt{\pi/3}}{\sqrt{a}}\operatorname{erf}(\sqrt{3a}).
\tag{56}
\end{align}

Now it is easy to get recurrence relation for $J_{n}$ in (53); we take $t=\sqrt{a k}\,x \Rightarrow dx=dt/(\sqrt{a k})$ and we come up to $\operatorname{erf}(\sqrt{a k})$ which gives the following general formula:
\begin{equation}
J_{n}=2+\frac{\sqrt{\pi}}{a}
\sum_{k=1}^{n}
\frac{(-1)^{k}}{\sqrt{k}}\binom{n}{k}\operatorname{erf}(\sqrt{a k}).
\tag{57}
\end{equation}

We produce in Table 1 a short table for comparison; we reused for this problem our approach with the same Pad\'e approximants and obtained the following as approximations:
\begin{equation}
I_{n}=
\frac{2}{2n+1}\left(\frac{4}{\pi}\right)^{n}
\,{}_2F_{1}\!\left(2n,2n+\tfrac{1}{2};2n+\tfrac{3}{2};-\tfrac{1}{3}\right),
\tag{58}
\end{equation}
\begin{equation}
I_{n}=
\frac{2}{2n+1}\left(\frac{4}{\pi}\right)^{n}
\,F_{1}\!\left(2n+\tfrac{1}{2};-2n,2n;2n+\tfrac{3}{2};\tfrac{1}{30},-\tfrac{3}{10}\right),
\tag{59}
\end{equation}
where $F_{1}$ is the Appell hypergeometric function of two variables.

\section{Series Representation of $T(x)$ Function over $[-1,1]$ Using Error Function Approximation and Pad\'e Approximant}

Recall
\begin{equation}
I_{k}=\int_{-1}^{1} x^{2k}\bigl(\operatorname{erf}(x)\bigr)^{k}\,dx
\tag{60}
\end{equation}
is $0$ if $k$ is odd. Thus, we need to focus on
\begin{equation}
I_{2k}=\int_{-1}^{1} x^{4k}\bigl(\operatorname{erf}(x)\bigr)^{2k}\,dx
\tag{61}
\end{equation}
which could be approximated, as we showed above in Section 3 to get (57), using
\begin{equation}
\bigl(\operatorname{erf}(x)\bigr)^{2}\approx 1-e^{-a x^{2}}
\quad\text{with}\quad
a=(1+\pi)^{2/3}\log^{2}(2),
\tag{62}
\end{equation}
making
\begin{equation}
I_{2k}=\int_{-1}^{1} x^{4k}\bigl(1-e^{-a x^{2}}\bigr)^{k}\,dx
\tag{63}
\end{equation}
to be developed using the binomial expansion. Therefore, in practice, we face the problem of
\begin{equation}
J_{n,k}=\int_{-1}^{1} x^{4k} e^{-n a x^{2}}\,dx
\tag{64}
\end{equation}
and the antiderivative
\begin{equation}
\int x^{4k} e^{-n a x^{2}}\,dx
=-\frac{1}{2}x^{4k+1}E_{1/2-2k}(a n x^{2}),
\tag{65}
\end{equation}
where the exponential integral function appears. Using the bounds, this reduces to
\begin{equation}
J_{n,k}=-E_{1/2-2k}(a n)
\tag{66}
\end{equation}
and leads to ``reasonable'' approximation as shown in Table 2.

Another approximation could be obtained using the simplest Pad\'e approximant \cite{KhaniShahmorad2012} of the error function
\begin{equation}
\operatorname{erf}(x)=\frac{2x}{\sqrt{\pi}\,(1+x^{2}/3)}
\tag{67}
\end{equation}
which would lead to
\begin{equation}
I_{2k}=\int_{-1}^{1} x^{4k}\bigl(\operatorname{erf}(x)\bigr)^{2k}\,dx
=\frac{2}{6k+1}\left(\frac{4}{\pi}\right)^{k}
\,{}_2F_{1}\!\left(2k,6k+\tfrac{1}{2};6k+\tfrac{3}{2};-\tfrac{1}{3}\right),
\tag{68}
\end{equation}
slightly less accurate than the previous one. Continuing with Pad\'e approximant
\begin{equation}
\operatorname{erf}(x)=\frac{2x}{\sqrt{\pi}}-\frac{x^{3}}{15\sqrt{\pi}}
\left(1+\frac{3x^{2}}{10}\right)^{-1}
\tag{69}
\end{equation}
we should get
\begin{equation}
I_{2k}=\int_{-1}^{1} x^{4k}\bigl(\operatorname{erf}(x)\bigr)^{2k}\,dx
=\frac{2}{6k+1}\left(\frac{4}{\pi}\right)^{k}
\,F_{1}\!\left(6k+\tfrac{1}{2};-2k,2k;6k+\tfrac{3}{2};\tfrac{1}{30},-\tfrac{3}{10}\right),
\tag{70}
\end{equation}
where the Appell hypergeometric function of two variables appears. Finally we conclude the series representation as follows:
\begin{equation}
I(t)=\int_{-1}^{1} \exp\bigl(-x^{2}\operatorname{erf}(x)\bigr)\,dx
\sim\sum_{k=0}^{\infty}\frac{(-1)^{k}}{k!}I_{2k}.
\tag{71}
\end{equation}

\section{Approximation of $T(x)$ Function by Means of a Polynomial}

\begin{lemma}
The function $f$ which is defined as
\begin{equation}
f(x)=T\!\left(\frac{b+a}{2}+\frac{b-a}{2}x\right),\qquad -1\le x\le 1
\tag{72}
\end{equation}
could be approximated by means of Chebytchev polynomial.
\end{lemma}

\paragraph{Proof.}
We may approximate the function $f$ on the interval $[-1,1]$ by using Ch\'ebyshev polynomials \cite{SolokloFarsangi2013} of the first kind. To this end, we choose some positive integer $n$ and we define the coefficients $c_{n}$ by the formula
\begin{equation}
c_{j}=\frac{2}{\pi}\int_{-1}^{1}\frac{T_{j}(x)}{\sqrt{1-x^{2}}}f(x)\,dx,
\qquad j=0,1,\dots,n.
\tag{73}
\end{equation}
Then the polynomial
\begin{equation}
P_{n}(x)=\frac{1}{2}c_{0}+\sum_{j=1}^{n}c_{j}T_{j}(x)
\tag{74}
\end{equation}
approximates $f(x)$ in the best possible way. Since
\begin{equation}
T(x)=f\!\left(\frac{a+b-2x}{a-b}\right)\quad\text{for }a\le x\le b
\tag{75}
\end{equation}
we see that the polynomial
\[
Q_{n}(a,b,x)=P_{n}\!\left(\frac{a+b-2x}{a-b}\right)
\]
is an approximant to $T(x)$ function on $[a,b]$. Calculations give
\begin{align}
Q_{11}\Bigl(0,\frac{3}{2},x\Bigr)
&=0.0137936039435\,x^{11}
-0.135129528505\,x^{10}
+0.548169602543\,x^{9}
-1.16161653976\,x^{8}
+1.31691631085\,x^{7}\nonumber\\
&\quad-0.746480407376\,x^{6}
+0.338453415662\,x^{5}
-0.370071852413\,x^{4}
+0.0133517048763\,x^{3}\nonumber\\
&\quad-0.00104123958376\,x^{2}
+1.00003172454\,x,
\tag{76}
\end{align}
and
\begin{align}
Q_{11}\Bigl(\frac{3}{2},3,x\Bigr)
&=-0.0000675632422240\,x^{11}
+0.00188305739843\,x^{10}
-0.0239397852528\,x^{9}
+0.183255163671\,x^{8}\nonumber\\
&\quad-0.937675010268\,x^{7}
+3.35913844398\,x^{6}
-8.55140470408\,x^{5}
+15.3046428836\,x^{4}\nonumber\\
&\quad-18.4622672665\,x^{3}
+13.5920479951\,x^{2}
-4.69093970289\,x
+1.04191571066.
\tag{77}
\end{align}
For both approximations the error is less than $10^{-6}$. Indeed, numerical integration gives
\begin{equation}
\bigl\|T(x)-Q_{11}(0,\tfrac{3}{2},x)\bigr\|
=\sqrt{\int_{0}^{3/2}\bigl(T(x)-Q_{11}(0,\tfrac{3}{2},x)\bigr)^{2}\,dx}
\approx 2.26\times 10^{-7}
\tag{78}
\end{equation}
and
\begin{equation}
\bigl\|T(x)-Q_{11}(\tfrac{3}{2},3,x)\bigr\|
=\sqrt{\int_{3/2}^{3}\bigl(T(x)-Q_{11}(\tfrac{3}{2},3,x)\bigr)^{2}\,dx}
\approx 3.66\times 10^{-10}.
\tag{79}
\end{equation}
Thus, we may evaluate the $T(x)$ function with high accuracy on the interval $[0,3]$. For $x>3$ we may use the following approximation formula in terms of the error function:
\begin{equation}
T(x)\approx \varphi(x)\stackrel{\mathrm{def}}{=}
\int_{0}^{3}\exp\bigl(-t^{2}\operatorname{erf}(t)\bigr)\,dt
+\frac{\sqrt{\pi}}{2}\bigl(\operatorname{erf}(x)-\operatorname{erf}(3)\bigr),
\qquad x\ge 3.
\tag{80}
\end{equation}
The quadratic mean error on $[3,100]$ is
\begin{equation}
\|T(x)-\varphi(x)\|
=\sqrt{\int_{3}^{100}\bigl(T(x)-\varphi(x)\bigr)^{2}\,dx}
\approx 2.02\times 10^{-8}.
\tag{81}
\end{equation}
Now we are ready to present application of $T(x)$ in probability and thermodynamics using one of the most important distributions which is called Maxwell–Boltzmann distribution.

\section{Application of $T(x)$ in Probability}

Let
\begin{equation}
F_{\lambda,\mu}(x)=\int_{0}^{x} e^{-\xi^{2}(\lambda+\mu\operatorname{erf}(\xi))}\,d\xi
\qquad (\lambda>0).
\tag{82}
\end{equation}
Define
\begin{equation}
c=\int_{0}^{\infty} e^{-\xi^{2}(\lambda+\mu\operatorname{erf}(\xi))}\,d\xi
\tag{83}
\end{equation}
and let
\begin{equation}
T_{\lambda,\mu}(x)\stackrel{\mathrm{def}}{=}c^{-1}F_{\lambda,\mu}(x),
\qquad x\ge 0.
\tag{84}
\end{equation}
The function $T(x)$ is the new special function we have studied in this paper. This function defines a cumulative probability distribution function (CDF) with probability distribution function (PDF)
\begin{equation}
f_{\lambda,\mu}(x)=e^{-x^{2}(\lambda+\mu\operatorname{erf}(x))},
\qquad x\ge 0.
\tag{85}
\end{equation}
Indeed, we have
\begin{equation}
T'_{\lambda,\mu}(x)=e^{-x^{2}(\lambda+\mu\operatorname{erf}(x))}>0,\qquad
T_{\lambda,\mu}(+\infty)=1.
\tag{86}
\end{equation}
The ODE for this function not involving the error function $\operatorname{erf}$ may be obtained by differentiating (84) twice and eliminating the expression containing that error function. This gives us the following ODE:
\begin{equation}
\sqrt{\pi}\,x\,y''(x)
=2e^{-x^{2}}y'(x)\bigl(\sqrt{\pi}e^{x^{2}}\log(y'(x))-\mu x^{3}\bigr).
\tag{87}
\end{equation}
Letting $\mu=0$ gives the ODE
\begin{equation}
x\,y''(x)=2y'(x)\log\bigl(y'(x)\bigr)
\tag{88}
\end{equation}
whose general solution is as follows:
\begin{equation}
y(x)=\frac{1}{2}e^{-c_{1}/2}\sqrt{\pi}\,\operatorname{erfi}\bigl(e^{c_{1}/2}x\bigr)+c_{2}.
\tag{89}
\end{equation}
If we compare (89) with
\begin{equation}
cF_{\lambda,0}(x)
=c\int_{0}^{x}e^{-\lambda\xi^{2}}\,d\xi
=\left(\int_{0}^{\infty}e^{-\lambda\xi^{2}}\,d\xi\right)^{-1}
\int_{0}^{x}e^{-\lambda\xi^{2}}\,d\xi
=\operatorname{erf}(\sqrt{\lambda}x)
\tag{90}
\end{equation}
we must have
\begin{equation}
\frac{1}{2}e^{-c_{1}/2}\sqrt{\pi}\,\operatorname{erfi}\bigl(e^{c_{1}/2}x\bigr)+c_{2}
=\operatorname{erf}(\sqrt{\lambda}x),
\tag{91}
\end{equation}
so that
\begin{equation}
c_{1}=i\pi+\log\Bigl(\frac{\pi}{4}\Bigr),\qquad c_{2}=0.
\tag{92}
\end{equation}
We showed that in the case when $\mu=0$ our function $T_{\lambda,0}(x)$ coincides with the error function $\operatorname{erf}(\sqrt{\lambda}x)$ with the value $\lambda=\pi/4$. When $\mu\ne 0$ we cannot obtain the solution to the ODE (87) in closed form. We may try a numerical procedure or another method to solve it. Our aim is to show how we may apply the new special function $T_{\lambda,\mu}(x)$ in probability and physics.

\subsection*{ Example}

We look for $\lambda$ and $\mu$ in order to adjust the error function by means of the function $y(x)=T_{\lambda,\mu}(x)$. To this end, we impose the following conditions:
\begin{equation}
\operatorname{erf}(1)=T_{\lambda,\mu}(1),\qquad
\operatorname{erf}'(1)=T'_{\lambda,\mu}(1).
\tag{93}
\end{equation}
Solving this system gives
\begin{equation}
\lambda=0.1671645,\qquad \mu=0.8449657.
\tag{94}
\end{equation}
The function $T_{\lambda,\mu}(x)$ converts into
\begin{equation}
T_{\lambda,\mu}(x)
=1.05021\int_{0}^{x}
\exp\bigl(-\xi^{2}(0.167164+0.844966\,\operatorname{erf}(\xi))\bigr)\,d\xi.
\tag{95}
\end{equation}
Plotting the two functions gives following picture as shown in Figure~10.

\begin{figure}[H]
  \centering
  \includegraphics[width=0.6\textwidth]{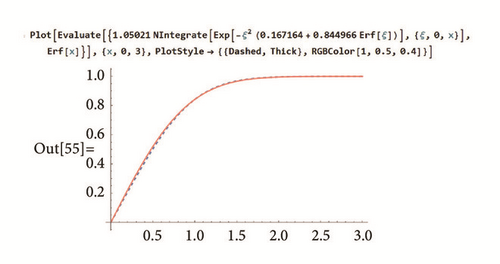}
  \caption{Adjusting the error function by means of the new special function.}
\end{figure}

\section{Application of $T(x)$ in Thermodynamics}

In physics (in particular in statistical mechanics), the Maxwell–Boltzmann distribution is a particular probability distribution named after James Clerk Maxwell and Ludwig Boltzmann. It was first defined and used for describing particle speeds in idealized gases, where the particles move freely inside a stationary container without interacting with one another, except for very brief collisions in which they exchange energy and momentum with each other or with their thermal environment. The term ``particle'' in this context refers to gaseous particles (atoms or molecules), and the system of particles is assumed to have reached thermodynamic equilibrium. The energies of such particles follow what is known as Maxwell–Boltzmann statistics, and the statistical distribution of speeds is derived by equating particle energies with kinetic energy. Mathematically, the Maxwell–Boltzmann distribution is the chi distribution with three degrees of freedom (the components of the velocity vector in Euclidean space), with a scale parameter measuring speeds in units proportional to the square root of $T/m$ (the ratio of temperature and particle mass); see Figure~11.

\begin{figure}[H]
  \centering
  \includegraphics[width=0.6\textwidth]{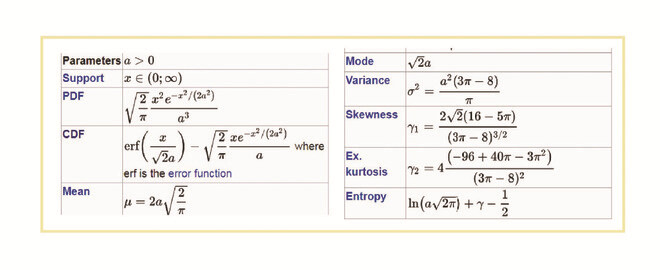}
  \caption{Maxwell–Boltzmann distribution.}
\end{figure}

The CDF for the Boltzmann distribution may be approximated by means of the new special function $T_{\lambda,\mu}(x)$ as follows:
\begin{equation}
\operatorname{erf}\!\left(\frac{x}{\sqrt{2}a}\right)
-\sqrt{\frac{2}{\pi}}\frac{x}{a}\exp\!\left(-\frac{x^{2}}{2a^{2}}\right)
\approx
T_{\lambda,\mu}(x)
-\sqrt{\frac{2}{\pi}}\frac{x}{a}\exp\!\left(-\frac{x^{2}}{2a^{2}}\right),
\tag{96}
\end{equation}
where $T_{\lambda,\mu}(x)$ is an approximation to $\operatorname{erf}(x/\sqrt{2}a)$ for some parameters $\lambda$ and $\mu$ depending on $a$. This approximation may be obtained in a similar way to what we illustrated in Example 1, Figure~11. On the other hand, in the case when $0<a\le 1$ we may approximate the CDF for the Maxwell–Boltzmann distribution as shown in Figure~12 for the value $a=0.75$.

\begin{figure}[H]
  \centering
  \includegraphics[width=0.6\textwidth]{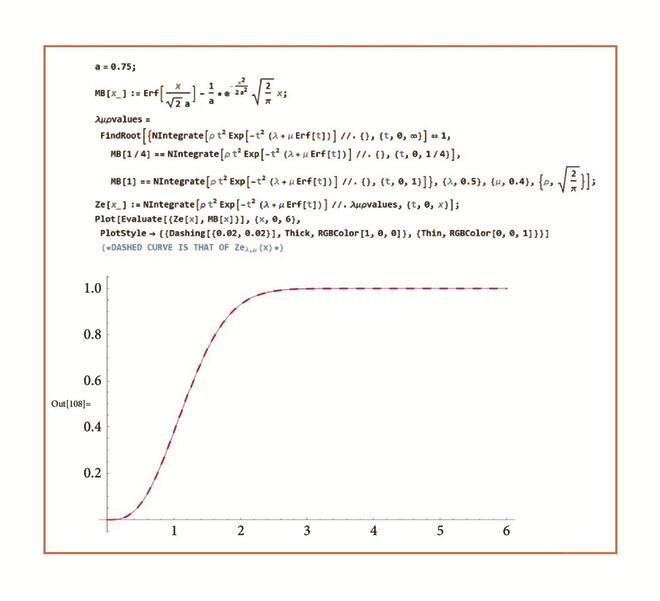}
  \caption{Approximation of the CDF for the Maxwell–Boltzmann distribution for $a=0.75$.}
\end{figure}

Finally this approximation by new special function showed that it may also be applied in thermodynamics to evaluate the average energy per particle in the circumstance where there is no energy-dependent density of states to skew the distribution, and the representation of probability for a given energy must be normalized to a probability of $1$ which holds using our new special function with two parameters as shown in (86).

\section{Conclusion}

We have studied a new probability distribution defined on $[0,+\infty)$ and we gave series representations for $T(x)$ function using Pad\'e approximant. Really we approximated the CDF for that distribution by means of Ch\'ebyshev polynomials and the error function. The methods we applied are suitable for approximating other CDF for probability distributions, since their CDF are bounded and they take values from $0$ to $1$. And it is well known that Ch\'ebyshev polynomials are the optimal ones for approximating continuous functions. On the other hand, it is also possible to approximate such functions by means of rational Ch\'ebyshev approximants. This technique may be used in future works.

\section*{Data Availability}

The data supporting this new special function is the Table of integral of error functions which is cited as one of the important reference in this paper and that function does not mentioned in any standard references related to the integral of error function.

\section*{Disclosure}

Zeraoulia Rafik present address is Department of Mathematics, High school-Timgad, Batna, Algeria.

\section*{Conflicts of Interest}

There are no conflicts of interest regarding the publication of this paper.

\section*{Authors’ Contributions}

Zeraoulia Rafik, Alvaro H. Salas, and David L. Ocampo are equally contributing authors.

\section*{Acknowledgments}

Alvaro H. Salas, Universidad Nacional de Colombia, Colombia, has received the Grant/Award Number 283, \url{http://agenciadenoticias.unal.edu.co/detalle/article/profesor-de-la-un-gano-premio-scopus-en-el-area-de-matematicas.html}. We would like to express our deep gratitude to Professor David L. Ocampo for improving the quality of the paper. My great salutation to Professor Alvaro H. Salas and to my parents and all my great salutations to my wife and to my second heart my son Taha Abd-Aldjalil. We would also like to thank Yuriy S and Claude Leibovici from Stack Exchange Math and Prof. Carlo Beenakker from MathOverflow for their contributions to the paper.

\section*{Endnotes}

1. In mathematics, the error function (also called the Gauss error function) is a special function (nonelementary) of sigmoid shape that occurs in probability, statistics, and partial differential equations describing diffusion. It is defined as $\operatorname{Erf}(x)=(2/\sqrt{\pi})\int_{0}^{x}e^{-t^{2}}dt$. Of course, it is closely related to the normal CDF $\Phi(x)=P(N<x)=(1/\sqrt{2\pi})\int_{-\infty}^{x}e^{-t^{2}/2}dt$ (where $N\sim N(0,1)$ is a standard normal) by the expression $\operatorname{Erf}=2\Phi(x\sqrt{2})-1$.

2. Cumulative distribution function for the normal distribution. In probability theory and statistics, the cumulative distribution function (CDF, also cumulative density function) of a real-valued random variable $X$, or just distribution function of $X$, evaluated at $x$, is the probability that $X$ will take a value less than or equal to $x$. If we have a quantity $A$ that takes some value at random, the cumulative density function $F(x)$ gives the probability that $X$ is less than or equal to $x$; that is,
\[
F(x)=P(A\le x).
\]
In the case of a continuous distribution, it gives the area under the probability density function from minus infinity to $x$. Cumulative distribution functions are also used to specify the distribution of multivariate random variables.


\begin{thebibliography}{99}

\bibitem{TeugelsSundt2004}
J.~L. Teugels and B.~Sundt,
\newblock {\em Encyclopedia of Actuarial Science}, vol.~1,
\newblock Wiley, 2004.

\bibitem{MaiDuyTranCong2003}
N.~Mai-Duy and T.~Tran-Cong,
\newblock Approximation of function and its derivatives using radial basis function networks,
\newblock {\em Applied Mathematical Modelling}, 27(3):197--220, 2003.

\bibitem{SuzukiSuzuki2003}
T.~Suzuki and T.~Suzuki,
\newblock Numerical integration error method for zeros of analytic functions,
\newblock {\em Journal of Computational and Applied Mathematics}, 152(1--2):493--505, 2003.

\bibitem{NgGeller1969}
E.~W. Ng and M.~Geller,
\newblock A table of integrals of the error functions,
\newblock {\em Journal of Research of the National Bureau of Standards}, 73B:1--20, 1969.

\bibitem{Liouville1838}
J.~Liouville,
\newblock Suite du M\'emoire sur la classification des Transcendantes, et sur l’impossibilit\'e d’exprimer les racines de certaines \'equations en fonction finie explicite des coefficients,
\newblock {\em Journal de Math\'ematiques Pures et Appliqu\'ees}, 3:523--546, 1838.

\bibitem{arxiv1702.08438}
Available at \texttt{https://arxiv.org/pdf/1702.08438.pdf}.

\bibitem{KhaniShahmorad2012}
A.~Khani and S.~Shahmorad,
\newblock An operational approach with Pade approximant for the numerical solution of non-linear Fredholm integro-differential equations,
\newblock {\em Scientia Iranica}, 19(6):1691--1698, 2012.

\bibitem{SolokloFarsangi2013}
H.~N. Soloklo and M.~M. Farsangi,
\newblock Chebyshev rational functions approximation for model order reduction using harmony search,
\newblock {\em Scientia Iranica}, 20(3):771--777, 2013.

\end{thebibliography}
\end{document}